\documentclass[reqno]{amsproc}
\usepackage[latin1]{inputenc}
\usepackage{amssymb}
\usepackage{hyperref}
\usepackage{amsmath}
\usepackage{latexsym}
\usepackage{cite}
\usepackage{amssymb}
\usepackage{amsfonts}
\usepackage{amscd}
\usepackage{xcolor}
\usepackage{amstext,amsmath,amssymb,amsfonts}
\newtheorem{theorem}{Theorem}

\newtheorem{corollary}[theorem]{Corollary}

\newtheorem{definition}[theorem]{Definition}
\newtheorem{example}[theorem]{Example}

\newtheorem{lemma}[theorem]{Lemma}

\newtheorem{proposition}[theorem]{Proposition}
\newtheorem{remark}[theorem]{Remark}

\newcommand{\K}{\mathbb {K}}

\newcommand{\A}{\mathcal{A}}

\newcommand{\beq}{\begin{eqnarray}}
\newcommand{\eeq}{\end{eqnarray}}
\newcommand{\beqs}{\begin{eqnarray*}}
\newcommand{\eeqs}{\end{eqnarray*}}
\newcommand{\bpro}{\begin{pro}}
\newcommand{\epro}{\end{pro}}
\newcommand{\blem}{\begin{lem}}
\newcommand{\elem}{\end{lem}}
\newcommand{\bdfn}{\begin{dfn}}
\newcommand{\edfn}{\end{dfn}}
\newcommand{\bcor}{\begin{cor}}
\newcommand{\ecor}{\end{cor}}
\newcommand{\bthm}{\begin{thm}}
\newcommand{\ethm}{\end{thm}}
\newcommand{\bex}{\begin{ex}}
\newcommand{\eex}{\end{ex}}
\newcommand{\brmk}{\begin{rmk}}
\newcommand{\ermk}{\end{rmk}}
\newcommand{\bpr}{\begin{pr}}
\newcommand{\epr}{\end{pr}}
\newcommand{\benum}{\begin{enumerate}}
\newcommand{\eenum}{\end{enumerate}}
\newcommand{\bitem}{\begin{itemize}}
\newcommand{\eitem}{\end{itemize}}

\chardef\bslash=`\\
\numberwithin{equation}{section}
\numberwithin{table}{section}
\numberwithin{theorem}{section}
\DeclareMathOperator{\id}{id}

\title[ 2-hom-associative bialgebras and hom-left symmetric dialgebras]{ 2-hom-associative bialgebras and hom-left symmetric dialgebras \footnote{Preprint: ICMPA-MPA/2015/08 } }

\author{Mahouton Norbert Hounkonnou$^\ast$}
\address[$\ast$]{University of Abomey-Calavi,
International Chair in Mathematical Physics and Applications,
ICMPA-UNESCO Chair, 072 BP 50, Cotonou, Rep. of Benin}
\email{norbert.hounkonnou@cipma.uac.bj, with copy to hounkonnou@yahoo.fr}
\author{Gb\^ev\`ewou Damien  Houndedji$^\dagger$}
\address[$\dagger$]{University of Abomey-Calavi,
International Chair in Mathematical Physics and Applications,
ICMPA-UNESCO Chair, 072 BP 50, Cotonou, Rep. of Benin}
\email{ houndedjid@gmail.com}
\begin{document}
\maketitle

\today

\bigskip
\begin{abstract}
From the definition and properties of  unital hom-associative algebras, and the use of  the Kaplansky's constructions, we develop new algebraic structures called {\it  2-hom-associative bialgebras, 2-hom-bialgebras, and 2-2-hom-bialgebras}. Besides, we devise a construction of   hom-left symmetric dialgebras and discuss  their main relevant properties.
\end{abstract}

\section{Introduction}
The generalizations of algebras became an exciting area of mathematics that was greatly developed over the past decades, motivated by their  interesting new algebraic identities and properties. 
A {\it $2$-associative algebra} is a linear space $V$ endowed with two associative products. This algebra  is said
 to be unital if  both the products have the same  unit $1.$ In particular, a {\it free 
 $2$-associative algebra over a vector space $V$} is a $2$-associative algebra, $2as(V),$ such that 
 any linear map from $V$ to a $2$-associative algebra $\A$ has a unique extension to a homomorphism 
 of $2$-associative algebras from $2as(V)$ to $\A$. See, for instance, 
 J-L. Loday and M. Ronco in \cite{[Loday3]}. A description of such an algebraic  structure in terms of 
 planar trees is given in \cite{[Loday3]}.  A free 2-associative algebra on one generator can be 
 identified with  non-commutative polynomials over the planar (rooted) trees. The operad of $2$-
 associative algebras is studied in \cite{[Loday3]}, and it turns out to be a Koszul operad.
 
 A bialgebra is a $\mathbb{K}$-vector space $V,$ where $\mathbb{K}$ is a field, equipped with an algebra structure
  given by a multiplication $\mu$ and a unit $\eta$, in addition to a coalgebra structure given by a
   comultiplication $\Delta$ and a counit $\epsilon$, such that there is a compatibility condition 
   between the two structures expressed by the fact that $\Delta$ and $\epsilon$ are algebra 
   morphisms, that is for all $x, y\in V:$
\beqs
\Delta(\mu(x \otimes y)) = \Delta(x) \bullet \Delta(y) \mbox{ and } \epsilon(\mu(x, y)) = \epsilon(x)\epsilon(y).
\eeqs

The multiplication $\bullet$ on $V\otimes V$ is the usual multiplication on tensor product as follows: 
\beqs
(x\otimes y)\bullet (x'\otimes y') = \mu(x\otimes x')\otimes \mu(y\otimes y').
\eeqs
The bialgebra is said infinitesimal if the compatibility condition is modified. The comultiplication 
is no more an algebra morphism. The condition is $\Delta\circ\mu= (\mu\otimes \id)\circ (\id\otimes 
\Delta) + (\id\otimes \mu)\circ (\Delta\otimes \id)$. This structure was introduced first by S. Joni 
and G.-C. Rota in \cite{[Joni-Rota]}. The basic theory was developed by M. Aguiar in \cite{[M1]} and in \cite{[M2]}. He also showed their intimate link with Rota-Baxter algebras, Loday's dendriform 
algebras, pre-Lie structure, and introduced the associative classical Yang-Baxter equation (see 
\cite{[M3]}).

 {\it A $2$-associative bialgebra $(V, \mu_{1}, \mu_{2}, \eta, \Delta, \epsilon)$} is a vector space $V$ 
equipped with two multiplications $\mu_{1}$ and $\mu_{2}$, a unit $\eta$, a comultiplication 
$\Delta,$ and a counit $\epsilon$ such that \cite{[Loday3]}
\begin{enumerate}
\item[(a)] $(V, \mu_{1}, \eta, \Delta, \epsilon)$ is a bialgebra,
\item[(b)] $(V, \mu_{2}, \eta, \Delta, \epsilon)$ is a unital infinitesimal bialgebra.
\end{enumerate}
In  \cite{[Dekkar-Makhlouf]}, the definitions of $2$-associative bialgebras, $2$-
bialgebras, and $2$-$2$-bialgebras  and their basic properties were explicitily given. In \cite{[Dekkar-Makhlouf]}, using the Kaplansky's constructions for 
a bialgebra, the authors gave a large class of these algebras.  They  classified, in 
dimensions $2$ and $3,$ the bialgebras leading to  the classification of infinitesimal 
bialgebra structures. They  also described the $2$ and $3$-dimensional $2$-associative 
bialgebras, $2$-bialgebras, and $2$-$2$-bialgebras. 

Further, J. L. Loday in \cite{[Loday1]} introduced a non-antisymmetric version of Lie algebras, called Leibniz algebra, whose bracket 
satisfies the Leibniz relation. The Leibniz 
rule, combined 
with the antisymmetry property, leads to the Jacobi identity. Therefore, the Lie algebras are anti-symmetric Leibniz 
algebras. In the same work, Loday also formulated an associative version of Leibniz algebras, 
called diassociative algebras, equipped with two bilinear and associative operations, which satisfy three axioms, all of them being various forms of the associative law. Recently \cite{[Felipe2]}, R. Felipe developed the concept of left-symmetric dialgebras, which includes, as a particular case,  the notion of dialgebras. This led to a new impulse to the construction of Leibniz algebras.

Hom-algebra structures first arose in quasi-deformations of Lie algebras of vector fields. Discrete 
modifications of vector fields via twisted derivations led to hom-Lie and quasi-hom-Lie structures, 
in which the Jacobi condition is twisted. Other interesting hom-type algebras of 
classical structures were studied as hom-associative algebras, hom-Lie admissible algebras\cite{[makhlouf]}, and, more 
generally, G-hom-associative algebras \cite{[makhlouf2]}, enveloping algebra of hom-Lie algebras\cite{[yau1]}, hom-Lie admissible hom-coalgebras and hom-Hopf algebras\cite{[makhlouf1]}, hom-alternative algebras, hom-Malcev algebras and hom-Jordan 
algebras \cite{[Donaldyau1]}, L-modules, L-comodules and hom-Lie quasi-bialgebras \cite{[Bakayoko1]}, and Laplacian of hom-Lie quasi-bialgebras \cite{[Bakayoko2]}.  

In this paper,  we devise a hom-type generalization of 2-associative algebras, 2-bialgebras, 2-associative bialgebras, 2-2-bialgebras and left symmetric 
dialgebras, leading to the concepts of  hom-bialgebras, 2-hom-associative algebras, 2-hom-bialgebras, 2-hom-
associative bialgebras, 2-2-hom-bialgebras and hom-left symmetric 
dialgebras,  respectively. Hom-type algebras are usually defined by twisting the defining axioms of a type of 
algebras by a certain twisting map. When the twisting map happens to be the identity map, we get 
an ordinary algebraic structure. Firstly, we introduce a hom-counital condition given as follows:
\beqs
(\varepsilon\otimes \alpha)\Delta(x)=(\alpha\otimes \varepsilon)\Delta(x)=\alpha^{2}(x), \forall x\in V.
\eeqs
This new hom-counital condition leads  to new definitions of counital hom-coassociative coalgebra and unital hom-bialgebra structures. Then, we introduce an unital infinitesimal hom-bialgebra condition given as follows:
  \beqs
  \Delta\circ\mu= (\mu\otimes \alpha)\circ(\alpha\otimes \Delta) + (\alpha\otimes \mu)\circ(\Delta\otimes \alpha) - \alpha^{2}\otimes \alpha.
  \eeqs   
This unital infinitesimal twisted condition lpermits  to define the unital infinitesimal hom-bialgebra structures. 
Then, we introduce the concepts of 2-hom-associative bialgebras, 2-hom-bialgebras and 2-2-hom 
bialgebras. Besides, we give a hom-algebra version of Kaplansky's constructions of hom-bialgebras in order to 
build unital $2$-associative hom-bialgebras, unital $2$-hom-bialgebras and unital $2$-$2$-hom-bialgebras. Finally, we define the notion of hom-left symmetric dialgebras generalizing the classical left symmetric dialgebras, and discuss their relevant properties.  

The paper is organized as follows. In section 2, we give the definitions of hom-bialgebra, 2-hom-associative algebra, 2-hom-associative bialgebra, 2-hom-bialgebra, 2-2-hom-bialgebra, and discuss their relevant properties. In section 3, we provide a hom-algebra version of Kaplansky's constructions of hom-bialgebras from a unital hom-associative algebra. We show that these constructions induce a large class of 2-hom-bialgebras, 2-hom-associative bialgebras, and 2-2-hom-bialgebras. In section 4, we introduce the notion of hom-left symmetric dialgebras, and discuss their relevant properties. 
\section{Unital 2-hom-associative bialgebras}
 In this section, we give the definitions of 2-hom-associative 
 algebras and unital infinitesimal hom-bialgebras. Then, we introduce the 
 concepts of 2-hom-associative bialgebras, 2-hom-bialgebras, and 2-2-hom bialgebras, 
 which require the notions of 2-hom-associative algebras, hom-bialgebras, and unital 
 infinitesimal hom-bialgebras.
 \subsection{Unital hom-associative algebra and counit hom-coassociative coalgebra}
 We briefly recall  the definitions of hom-associative algebra, unital hom-associative algebra, and hom-coassociative coalgebra and their properties pertaining to this work. Then, we provide a  hom-counital condition  allowing to  define a counital hom-coassociative coalgebra. 
 \begin{definition}\cite{[makhlouf2]}
A hom-associative algebra is a triple $(V, \mu, \alpha)$ consisting of  a linear space V, a bilinear map $\mu: V\times V\rightarrow V$ and a homomorphism $\alpha: V\rightarrow V$ satisfying
\beq
\alpha\circ\mu= \mu\circ\alpha^{\otimes^{2}}:=  \mu\circ(\alpha\otimes\alpha)\mbox{(multiplicativity)}
\eeq
\beq
\mu\circ(\alpha\otimes \mu)= \mu\circ(\mu\otimes\alpha)\mbox{(hom-associativity)}.
\eeq
\end{definition}
\begin{example}\cite{[makhlouf2]}
Let $\A$ be a $3$-dimensional vector space over $\K$ with a basis $\lbrace e_{1}, e_{2}, e_{3}\rbrace$. The following multiplication $\mu$ and linear map $\alpha$ on $\A$ define a hom-associative algebra:
\beqs
&&\mu(e_{1}, e_{1})= ae_{1};\ \  \ \  \mu(e_{2}, e_{2})= ae_{2}\cr
&&\mu(e_{1}, e_{2})=\mu(e_{2}, e_{1})= ae_{2}; \mu(e_{2}, e_{3})=be_{3}\cr
&&\mu(e_{1}, e_{3})=\mu(e_{3}, e_{1})= be_{3},
\eeqs
where $a, b$ are parameters in $\K,$ and
\beqs
\alpha(e_{1})=ae_{1};\  \ \alpha(e_{2})=ae_{2};\ \ \alpha(e_{3})=be_{3}.
\eeqs
When $a\neq b$ and $b\neq 0$, the equality $\mu(\mu(e_{1}, e_{1}), e_{3})- \mu(e_{1}, \mu(e_{1}, e_{3}))= (a-b)be_{3}$ makes this algebra to be non-associative.
\end{example}
\begin{definition}\cite{[makhlouf2]}
Let $(V, \mu, \alpha)$ and $(V', \mu', \alpha')$ be two hom-associative algebras. A liner map $f: V\rightarrow V'$ is a morphism of hom-associative algebras if
\beq
\mu' \circ(f\otimes f)= f\circ\mu \mbox{ and } f\circ\alpha= \alpha' \circ f.
\eeq
In particular, the hom-associative algebras $(V, \mu, \alpha)$ and $(V', \mu', \alpha')$ are isomorphic if $f$ is a bijective linear map such that
\beq
\mu= f^{-1}\circ\mu'\circ(f\otimes f)  \alpha= f^{-1}\circ\alpha'\circ f.
\eeq 
\end{definition}
\begin{theorem}\label{t1}
Let $(V, \mu)$ be an associative algebra and $\alpha: V\rightarrow V$ be associative algebra
 endomorphism. Then $V_{\alpha}=(V, \alpha\circ\mu, \alpha),$ is a hom-associative algebra.
  Moreover, suppose that $(V', \mu')$ is another associative algebra and $\alpha': V'\rightarrow
   V'$ is an associative algebra endomorphism. If $f: V\rightarrow V'$ is an associative algebra
    morphism that satisfies $f\circ\alpha= \alpha'\circ f,$ then 
\beqs
f: (V, \alpha\circ \mu, \alpha)\rightarrow (V', \alpha'\circ \mu', \alpha')
\eeqs
is a morphism of hom-associative algebras.
\end{theorem}
\begin{definition}\cite{[makhlouf2]}
The hom-associative algebra is called unital if it admits a unity, ie. an element $e \in V$ such that $\mu(x, e)=\mu(e, x)=x$ for all $x\in V.$ The unity may also be expressed by a linear map $\eta: \K\rightarrow V$ defined by $\eta(c)=ce,$ for all $c\in \K$. Then, we denote the unital hom-associative algebra by $(V, \mu, \eta, \alpha).$
\end{definition}
Let $(V, \mu, \eta, \alpha)$ and $(V', \mu', \eta', \alpha')$ be two unital hom-associative
 algebras. Then the morphisms $f$ of unital hom-associative algebras are also required  to preserve the
  unital structures, i.e, to satisfy $f\circ\eta= \eta',$ and in the definition of isomorphism of
   unital hom-associative algebras, it is also required that $\eta= f^{-1}\circ\eta'.$
   \begin{example}
Let $\A$ be a $2$-dimensional vector space over $\K$ with a basis $\lbrace e_{1}, e_{2}\rbrace$. The following multiplication $\mu$ and linear map $\alpha$ on $\A$ define a unital hom-associative algebra, i. e.
\beqs
&&\mu(e_{1}, e_{1})= e_{1};  \ \  \mu(e_{2}, e_{2})= \mu(e_{1}, e_{2})=\mu(e_{2}, e_{1})= e_{2} \mbox{ and }\cr
&&\alpha(e_{1})=\alpha(e_{2})=e_{2}.
\eeqs
\end{example}
\begin{remark}
The tensor product of hom-associative algebras $(V, \mu, \alpha, \eta)$ and $(V', \mu', \alpha', \eta')$ is defined in an obvious way as the hom-associative algebra $(V\otimes V', \mu\otimes\mu', \alpha\otimes\alpha').$ 
\end{remark}
 \begin{definition}\cite{[yau3]}
A hom-coassociative coalgebra is a triple $(V, \Delta, \alpha)$ consisting of a linear space V, a linear map $\Delta: V\rightarrow V\otimes V,$ and a homomorphism $\alpha: V\rightarrow V$ satisfying
\beq
\alpha^{\otimes^{2}}\circ\Delta =\Delta\circ \alpha \mbox{(comultiplicativity)}
\eeq
\beq
(\alpha\otimes \Delta)\circ \Delta= (\Delta\otimes\alpha)\circ\Delta\mbox{(hom-coassociativity)}.
\eeq
\end{definition}
   \begin{example}
Let $\A$ be a $2$-dimensional vector space over $\K$ with a basis $\lbrace e_{1}, e_{2}\rbrace$. The following coproduct $\Delta$  and linear map $\alpha$ on $\A$ define a hom-coassociative coalgebra, i. e.
\beqs
&&\Delta(e_{1})=\Delta(e_{2})= e_{2}\otimes e_{2} \mbox{ and } \alpha(e_{1})=\alpha(e_{2})=e_{2}.
\eeqs
\end{example}
Now, we are in a position to give a new definition of a counital hom-coassociative coalgebra.  
Recalling that a counital coassociative coalgebra is a triple $(V, \Delta, \varepsilon)$ which is a coalgebra satisfying the counital condition:
\beqs
(\varepsilon\otimes id_{V})\Delta(x)=(id_{V}\otimes \varepsilon)\Delta(x)=id_{V}(x), \forall x\in V,
\eeqs
we get the following statement:
\begin{definition}
A counital hom-coassociative coalgebra is a quadruple $(V, \Delta, \varepsilon, \alpha)$ such that the triple $(V, \Delta, \alpha)$ is a hom-coassociative coalgebra satisfying the hom-counital condition
\beq\label{hom-counit condition}
(\varepsilon\otimes \alpha)\Delta(x)=(\alpha\otimes \varepsilon)\Delta(x)=\alpha^{2}(x), \forall x\in V.
\eeq
\end{definition} 
\begin{remark}
The condition (\ref{hom-counit condition}) can also  be replaced by
\beq
(\varepsilon\otimes id_{V})\Delta(x)=(id_{V}\otimes \varepsilon)\Delta(x)=\alpha(x), \forall x\in V.
\eeq
\end{remark}
   \begin{example}
Let $\A$ be a $2$-dimensional vector space over $\K$ with a basis $\lbrace e_{1}, e_{2}\rbrace$. The following coproduct $\Delta$  and linear maps $\alpha, \varepsilon$ on $\A$  define a counital  hom-coassociative coalgebra,  i. e. 
\beqs
&&\Delta(e_{1})=\Delta(e_{2})= e_{2}\otimes e_{2};\ \alpha(e_{1})=\alpha(e_{2})=e_{2}  \mbox{ and } \varepsilon(e_{1})=0, \varepsilon(e_{2})= 1.
\eeqs
\end{example}
 \subsection{Unital hom-bialgebra and unital infinitesimal hom-bialgebra}
 In the sequel, we define the unital hom-bialgebra and unital infinitesimal hom-bialgebra, where the
 comultiplication is a hom-algebra morphism, and the unital infinitesimal hom-bialgebra condition is modified.
 
  \begin{definition}
  A unital hom-bialgebra is a system $(V, \mu, \eta, \Delta, \varepsilon, \alpha)$, where   
  $\mu: V\otimes V\rightarrow V$(multiplication), $\eta: \K\rightarrow V$(unit), $\Delta: 
  V\rightarrow V\otimes V$ (comultiplication), $\varepsilon : V\rightarrow \K$ (counit), and $\alpha: 
  V\rightarrow V$ (endomorphism) are linear maps satisfying:
  \begin{enumerate}
  \item the quadruple $(V, \mu, \eta, \alpha)$ is a unital hom-associative algebra;
  \item the quadriple $(V, \Delta, \varepsilon, \alpha)$ is a counital hom-coassociative coalgebra;
  \item the compatibility condition is expressed by the following three identities:
  \begin{enumerate}
  \item $\Delta(\mu(x\otimes y))= \Delta(x)\bullet\Delta(y), \forall x, y \in V$,
  \item $\alpha^{\otimes^{2}}\circ\Delta =\Delta\circ \alpha,$
  \item $\varepsilon(\mu(x\otimes y))= \varepsilon(x)\varepsilon(y),$
  \item $\varepsilon\circ\alpha(x)=\varepsilon(x)$.
  \end{enumerate}
\end{enumerate}   
  \end{definition}
  \begin{remark}
  The conditions (3.a) and (3.b) mean that $\Delta$ is a homomorphism of the hom-associative algebras $(V, \mu, \alpha)$ and $(V\otimes V, \bullet, \alpha\otimes\alpha),$ while the conditions (3.c) and (3.d) tell us that $\varepsilon$ is a homomorphism of the hom-associative algebras $(V, \mu, \alpha)$ and $(\mathbb{K}, \cdot, id_{\mathbb{K}})$.  
  \end{remark}
     \begin{example}\label{ex1}
Let $\A$ be a $2$-dimensional vector space over $\K$ with a basis $\lbrace e_{1}, e_{2}\rbrace$. The product $\mu$, the coproduct $\Delta$ and the linear maps $\alpha, \varepsilon$ on $\A$ define a unital hom-bialgebra, i. e. 
\beqs
&&\mu(e_{1}, e_{1})= e_{1};  \ \  \mu(e_{2}, e_{2})= \mu(e_{1}, e_{2})=\mu(e_{2}, e_{1})= e_{2};\cr
&&\Delta(e_{1})=\Delta(e_{2})= e_{1}\otimes e_{1}; \alpha(e_{1})=\alpha(e_{2})=e_{2} \mbox{ and } \varepsilon(e_{1})= 1, \varepsilon(e_{2})= 0.
\eeqs
\end{example}
 J-L. Loday and M. Ronco   defined a {\it unital infinitesimal bialgebra} $(V, \mu, \eta, \Delta, \varepsilon)$ as a vector space $V$ equipped with a unital associative multiplication $\mu$ and a counital coassociative comultiplication $\Delta$ which are related by the unital infinitesimal relation:
  \beqs
  \Delta\circ\mu= (\mu\otimes id_{V})\circ(id_{V}\otimes \Delta) + (id_{V}\otimes \mu)\circ(\Delta\otimes id_{V}) - id_{V}\otimes id_{V}.
  \eeqs
Here, we formulate  the definition of a unital infinitesimal hom-bialgebra as follows:
\begin{definition}
A unital infinitesimal hom-bialgebra $(V,\mu, \eta, \Delta, \varepsilon, \alpha)$ is a $\mathbb{K}$-vector space $V$ equipped with a unital hom-associative multiplication $\mu$ and a counital hom-coassociative comultiplication $\Delta$ which are related by the unital hom-infinitesimal relation 
  \beq
  \Delta\circ\mu= (\mu\otimes \alpha)\circ(\alpha\otimes \Delta) + (\alpha\otimes \mu)\circ(\Delta\otimes \alpha) - \alpha^{2}\otimes \alpha.
  \eeq
\end{definition}
\begin{example}
The unital hom-bialgebra, given in Example \ref{ex1}, is a unital infinitesimal hom-bialgebra. 
\end{example}
 \subsection{2-hom-associative algebra}
The $2$-hom-associative algebras generalize the $2$-associative algebras in the sense  where the associativity laws are twisted.
\begin{definition}
A $2$-hom-associative algebra over $\K$ is a vector space equipped with two hom-associative structures. A $2$-hom-associative algebra is said to be unital if there is a unit $1$ which is a unit for both operations.
\end{definition}
   \begin{example}
Let $\A$ be a $2$-dimensional vector space over $\K$ with a basis $\lbrace e_{1}, e_{2}\rbrace$. The following multiplications $\mu_{1}, \mu_{2}$ and linear map $\alpha$ on $\A$ define a unital $2$- hom-associative algebra:
\beqs
&&\mu_{1}(e_{1}, e_{1})= e_{1};  \ \  \mu_{1}(e_{2}, e_{2})= \mu_{1}(e_{1}, e_{2})=\mu_{1}(e_{2}, e_{1})= e_{2};\cr
&& \mu_{2}(e_{1}, e_{1})= e_{1};  \ \  \mu_{2}(e_{2}, e_{2})= 0,  \mu_{2}(e_{1}, e_{2})=\mu_{2}(e_{2}, e_{1})= e_{2} \mbox{ and }\cr
&&\alpha(e_{1})=\alpha(e_{2})=e_{2}.
\eeqs
\end{example}
\begin{definition}
Let $(V, \mu_{1}, \mu_{2}, \alpha)$ and $(V', \mu'_{1}, \mu'_{2}, \alpha')$ be two $2$-hom-associative algebras. A liner map $f: V\rightarrow V'$ is a morphism of $2$-hom-associative algebras if
\beqs
\mu'_{1} \circ(f\otimes f)= f\circ\mu_{1},\ \ \mu'_{2} \circ(f\otimes f)= f\circ\mu_{2} \mbox{ and } f\circ\alpha= \alpha' \circ f.
\eeqs
In particular, the $2$-hom-associative algebras $(V, \mu_{1}, \mu_{2}, \alpha)$ and $(V', \mu'_{1}, \mu'_{2}, \alpha')$ are isomorphic if $f$ is a bijective linear map such that
\beqs
\mu_{1}= f^{-1}\circ\mu'_{1}\circ(f\otimes f),\ \ \mu_{2}= f^{-1}\circ\mu'_{2}\circ(f\otimes f) \mbox{ and } \alpha= f^{-1}\circ\alpha'\circ f.
\eeqs 
\end{definition}
\begin{theorem}
Let $(V, \mu_{1}, \mu_{2})$ be a $2$-associative algebra and $\alpha: V\rightarrow V$ be an associative algebra
 endomorphism. Then, $V_{\alpha}=(V, \alpha\circ\mu_{1}, \alpha\circ\mu_{2}, \alpha)$ is a $2$-hom-associative algebra.
  Moreover, suppose that $(V', \mu'_{1}, \mu'_{2})$ is another $2$-associative algebra and $\alpha': V'\rightarrow V'$ an associative algebra endomorphism. If $f: V\rightarrow V'$ is an associative algebra
    morphism that satisfies $f\circ\alpha= \alpha'\circ f,$ then 
\beqs
f: (V, \alpha\circ\mu_{1}, \alpha\circ\mu_{2}, \alpha)\rightarrow (V', \alpha\circ\mu'_{1}, \alpha\circ\mu'_{2}, \alpha')
\eeqs
is a morphism of $2$-hom-associative algebras.
\end{theorem} 
\textbf{Proof:}
It stems from Theorem\ref{t1}.

$ \hfill \square $

  \subsection{Unital 2-hom-associative bialgebra}
  We give the notion of unital $2$-hom-associative bialgebras generalizing unital $2$-associative bialgebras. Then, a unital $2$-hom-associative bialgebra is defined as follows:
\begin{definition}
A unital 2-hom-associative bialgebra $(V, \mu_{1}, \mu_{2}, \eta, \Delta, \varepsilon, 
\alpha)$ is a vector space $V$ equipped with two multiplications $\mu_{1}$ and $\mu_{2},$ 
a unit $\eta,$ a comultiplication $\Delta$, a counit $\varepsilon,$ and a linear map $\alpha: 
V\rightarrow V$ such that
\begin{enumerate}
\item $(V, \mu_{1}, \eta, \Delta, \varepsilon, \alpha)$ is a unital hom-bialgebra,
\item $(V, \mu_{2}, \eta, \Delta, \varepsilon, \alpha)$ is a unital infinitesimal hom-bialgebra.
\end{enumerate}
\end{definition}
     \begin{example}\label{ex2}
Let $\A$ be a $2$-dimensional vector space over $\K$ with a basis $\lbrace e_{1}, e_{2}\rbrace$. The product $\mu$, the coproduct $\Delta$ and the linear maps $\alpha$ and $\varepsilon$ given by
\beqs
&&\mu_{1}(e_{1}, e_{1})= e_{1};  \ \  \mu_{1}(e_{2}, e_{2})= \mu_{1}(e_{1}, e_{2})=\mu_{1}(e_{2}, e_{1})= e_{2};\cr
&& \mu_{2}(e_{1}, e_{1})= e_{1};  \ \  \mu_{2}(e_{2}, e_{2})= 0,  \mu_{2}(e_{1}, e_{2})=\mu_{2}(e_{2}, e_{1})= e_{2};\cr
&&\Delta(e_{1})=\Delta(e_{2})= e_{1}\otimes e_{1}; \alpha(e_{1})=\alpha(e_{2})=e_{2}; \varepsilon(e_{1})= 1, \varepsilon(e_{2})= 0
\eeqs
define a unital 2-hom-associative bialgebra structure on $\A.$
\end{example}

\begin{definition}
Let $(V, \mu_{1}, \mu_{2}, \eta, \Delta, \varepsilon, \alpha)$ and $(V', \mu'_{1}, \mu'_{2}, \eta', \Delta', \varepsilon', \alpha')$ be two unital\\ $2$-associative hom-bialgebras. A linear map $f: V\rightarrow V'$ is a morphism of unital $2$-associative hom-bialgebras if $
\mu'_{1}\circ(f\otimes f)= f\circ\mu_{1}, \mu'_{2}\circ(f\otimes f)= f\circ\mu_{2}, f\circ\eta= \eta', (f\otimes f)\circ \Delta=\Delta'\circ f, \varepsilon'=\varepsilon\circ f, f\circ\alpha=\alpha'\circ f.$
\end{definition}

 \subsection{Unital 2-hom-bialgebra}
 \begin{definition}
A unital 2-hom-bialgebra $(V, \mu_{1}, \mu_{2}, \eta, \Delta_{1}, \Delta_{2}, \varepsilon_{1}, \varepsilon_{2}, \alpha)$ is a vector space $V$ equipped with two multiplications $\mu_{1}, \mu_{2},$ the unit $\eta,$ two comultiplications $\Delta_{1}, \Delta_{2},$ two counits $\varepsilon_{1}, \varepsilon_{2},$ and a linear map $\alpha: V\rightarrow V$ such that:

 $(V, \mu_{1}, \eta, \Delta_{1}, \varepsilon_{1}, \alpha),$ $(V, \mu_{2}, \eta, \Delta_{2}, \varepsilon_{2}, \alpha),$ $(V, \mu_{1}, \eta, \Delta_{2}, \varepsilon_{2}, \alpha),$ and $(V, \mu_{2}, \eta, \Delta_{1}, \varepsilon_{1}, \alpha)$ 

are unital hom-bialgebras.
\end{definition}
\begin{example}
The unital 2-hom-associative bialgebra, given in Example \ref{ex2}, is a unital 2-hom-bialgebra, where $\Delta=\Delta_{1}=\Delta_{2}$.
\end{example}
The unital $2$-hom-bialgebra is called of type (1-1), (resp. of type (2-2)), if the two multiplications and the two comultiplications are identical, (resp. distinct).
The unital $2$-hom-bialgebra is called of type (1-2), (resp. of type (2-1)), if the two multiplications  are identical, (resp. distinct),  and the two comultiplications are distinct, (resp. identical).

\begin{proposition}
Let $(V,\mu, \eta, \Delta, \varepsilon, \alpha)$ be a unital hom-bialgebra, then $(V, \mu, \mu, \eta, \Delta, \Delta, \varepsilon, \alpha)$ and $(V, \mu, \mu^{op}, \eta, \Delta, \Delta^{cop}, \varepsilon, \alpha)$ are unital $2$-hom-bialgebras, where $\mu^{op}(x\otimes y)= \mu(y\otimes x)$ and $\Delta^{cop}(x)= \tau\circ \Delta(x),$ with $\tau(x\otimes y)= y\otimes x.$ The first unital $2$-hom-bialgebra is of type ($1$-$1$), and the second is of type ($2$-$2$). 
\end{proposition}
\textbf{Proof:} It comes from a direct computation.

$ \hfill \square $
 \subsection{Unital 2-2-hom-bialgebra}
\begin{definition}
A unital $2$-$2$-hom-bialgebra $(V, \mu_{1}, \mu_{2}, \eta, \Delta_{1}, \Delta_{2}, \varepsilon_{1}, \varepsilon_{2}, \alpha)$ is a vector space $V$ equipped with two multiplications $\mu_{1}, \mu_{2},$ two comultiplications $\Delta_{1}, \Delta_{2},$ two counits $\varepsilon_{1}, \varepsilon_{2}$, one unit $\eta,$ and a linear map $\alpha: V\rightarrow V$ such that 
\begin{enumerate}
\item $(V, \mu_{1}, \eta, \Delta_{1}, \varepsilon_{1}, \alpha)$ and $(V, \mu_{2}, \eta, \Delta_{2}, \varepsilon_{2}, \alpha)$ are unital hom-bialgebras,
\item $(V, \mu_{1}, \eta, \Delta_{2}, \varepsilon_{2}, \alpha)$ and $(V, \mu_{2}, \eta, \Delta_{1}, \varepsilon_{1}, \alpha)$ are unital infinitesimal hom-bialgebras.
\end{enumerate}
\end{definition}
\begin{example}
The unital 2-hom-associative bialgebra, given in Example \ref{ex2}, is a unital 2-2-hom-bialgebra, where $\Delta=\Delta_{1}=\Delta_{2}$.
\end{example}
A unital $2$-$2$-hom-bialgebra is called of type (1-1), (resp. of type (2-2)),  if the two multiplications and the two comultiplications are identical, (resp. distinct).
A unital $2$-$2$-hom-bialgebra is called of type (1-2), (resp. of type (2-1)),  if the two multiplications  are identical, (resp. distinct), and the two comultiplications are distinct, (resp. identical).

The definition of unital $2$-$2$-hom-bialgebra morphism is similar to that of unital $2$-hom-bialgebra morphism.
\section{Kaplansky's construction of hom-bialgebras}
In this section, we give a hom-algebra version of Kaplansky's constructions of hom-bialgebras in order to 
build unital $2$-associative hom-bialgebras, unital $2$-hom-bialgebras, and unital $2$-$2$-hom-bialgebras. The following statement is in order.
 \begin{proposition}\label{k1}
 Let $\mathcal{A}=(V, \mu, \eta, \alpha)$ be a unital hom-associative algebra,  where $e_{2}:= 
 \eta(1)$ is the unit. Let $\tilde{V}$ be the vector space spanned by $V$ and $e_{1}, 
 \tilde{V}=span(V, e_{1})$. We have the unital hom-bialgebra $\mathcal{K}_{1}
 (\mathcal{A}):=(\tilde{V}, \mu_{1}, \eta_{1}, \Delta_{1}, \varepsilon_{1}, \alpha_{1}),$ where the 
 multiplication $\mu_{1}$ is defined by:
 \beqs
&& \mu_{1}(e_{1}\otimes x)= \mu_{1}(x\otimes e_{1})= x \ \ \forall x\in \tilde{V},\cr
&& \mu_{1}(x\otimes y)= \mu(x\otimes y) \ \ \forall x, y\in V,
 \eeqs
 the unit $\eta_{1}$ is given by $\eta_{1}(1)= e_{1},$ while the comultiplication $\Delta_{1},$  the counit $\varepsilon_{1}, $  and the linear map $\alpha_{1}$ are defined by,  $\forall x\in V$:
 \beqs
 &&\Delta_{1}(e_{1})=e_{1}\otimes e_{1}\cr
 &&\Delta_{1}(x)= \alpha(x)\otimes e_{1} + e_{1}\otimes \alpha(x)- e_{2}\otimes \alpha(x) 
 \eeqs
 \beqs
 \varepsilon_{1}(e_{1})=1, \varepsilon_{1}(x)=0 
 \eeqs
 \beqs
 \alpha_{1}(e_{1})= e_{1}, \alpha_{1}(e_{2})= e_{2}, \alpha_{1}(x)=\alpha(x)
 \eeqs
respectively.
 \end{proposition}
 \textbf{Proof:}
 \begin{enumerate}
\item[$\diamond$] 
 $(V, \mu, \eta, \alpha)$ is a unital hom-associative algebra and 
 \begin{center}
 $\mu_{1}(\alpha_{1}(e_{1}), \mu_{1}(e_{1}, e_{1}))= e_{1}= \mu_{1}(\mu_{1}(e_{1}, e_{1}), \alpha_{1}(e_{1}))$.
 \end{center}
 Then $\forall x, y, z\in \tilde{V}$, $\mu_{1}(\alpha_{1}(x), \mu_{1}(y, z))= \mu_{1}(\mu_{1}(x, y), \alpha_{1}(z))$. Hence $(\tilde{V}, \mu_{1}, \eta_{1}, \alpha_{1})$ is a unital hom-associative algebra.
 \item[$\diamond$]
 We have
 \beqs
 (\alpha_{1}\otimes\Delta_{1})\circ\Delta_{1}(x)&=&(\alpha_{1}\otimes\Delta_{1})(\alpha_{1}(x)\otimes e_{1} + e_{1}\otimes\alpha_{1}(x) -e_{2}\otimes\alpha_{1}(x))\cr
 &=&\alpha^{2}_{1}(x)\otimes\Delta_{1}(e_{1}) + e_{1}\otimes\Delta_{1}(\alpha_{1}(x)) - \alpha^{2}_{1}(e_{2})\otimes\Delta_{1}(\alpha_{1}(x))\cr
 &=&\alpha^{2}_{1}(x)\otimes e_{1}\otimes e_{1} + e_{1}\otimes\alpha^{2}_{1}(x)\otimes e_{1} + e_{1}\otimes e_{1}\otimes\alpha^{2}_{1}(x) -e_{1}\otimes e_{2}\otimes\alpha^{2}_{1}(x)\cr
 && -e_{2}\otimes\alpha^{2}_{1}(x)\otimes e_{1}- e_{2}\otimes e_{1}\otimes\alpha^{2}_{1}(x) + e_{2}\otimes e_{2}\otimes\alpha^{2}_{1}(x)\cr
 &=&\Delta_{1}(\alpha_{1}(x))\otimes e_{1} + \Delta_{1}(e_{1})\otimes\alpha^{2}_{1}(x) + \Delta_{1}(e_{2})\otimes\alpha^{2}_{1}(x)\cr
 &=& (\Delta_{1}\otimes \alpha_{1})\circ\Delta_{1}(x), 
 \eeqs
and
\beqs
&&(\varepsilon_{1}\otimes\alpha_{1})\circ\Delta_{1}(x)=(\varepsilon_{1}\otimes\alpha_{1})(\alpha_{1}(x)\otimes e_{1} + e_{1}\otimes\alpha_{1}(x)- e_{2}\otimes\alpha_{1}(x))
=1\otimes\alpha^{2}_{1}(x)=\alpha^{2}_{1}(x) \mbox{ and }\cr
&&(\alpha_{1}\otimes \varepsilon_{1})\circ\Delta_{1}(x)=(\alpha_{1}\otimes \varepsilon_{1})(\alpha_{1}(x)\otimes e_{1} + e_{1}\otimes\alpha_{1}(x)- e_{2}\otimes\alpha_{1}(x))= \alpha^{2}_{1}(x)\otimes 1=\alpha^{2}_{1}(x).
\eeqs 
Then $(\varepsilon_{1}\otimes\alpha_{1})\circ\Delta_{1}(x)= (\alpha_{1}\otimes \varepsilon_{1})\circ\Delta_{1}(x)= \alpha^{2}_{1}(x).$ Hence, we can conclude that $(\tilde{V}, \Delta_{1}, \varepsilon_{1}, \alpha_{1})$ is a counital hom-coassociative coalgebra. 
 \item[$\diamond$]
\beqs
\Delta_{1}(x)\bullet\Delta_{1}(y)&=&(\alpha_{1}(x)\otimes e_{1} + e_{1}\otimes\alpha_{1}(x) -e_{2}\otimes\alpha_{1}(x))\bullet(\alpha_{1}(y)\otimes e_{1} + e_{1}\otimes\alpha_{1}(y) -e_{2}\otimes\alpha_{1}(y))\cr
&=&\mu(\alpha(x), \alpha(y))\otimes e_{1} +\alpha(x)\otimes\alpha(y) - \alpha(x)\otimes\alpha(y) + \alpha(y)\otimes\alpha(x) + e_{1}\otimes\mu(\alpha(x), \alpha(y))\cr
&& - e_{2}\otimes\mu(\alpha(x), \alpha(y)) - \alpha(y)\otimes\alpha(x) - e_{2}\otimes\mu(\alpha(x), \alpha(y)) +  e_{2}\otimes\mu(\alpha(x), \alpha(y))\cr
&=&\mu_{1}(\alpha(x), \alpha(y))\otimes e_{1} + e_{1}\otimes\mu_{1}(\alpha(x), \alpha(y)) - e_{2}\otimes\mu_{1}(\alpha(x), \alpha(y))\cr
&=&\mu_{1}\circ\alpha^{\otimes^{2}}_{1}(x\otimes y)\otimes e_{1} + e_{1}\otimes\mu_{1}\circ\alpha^{\otimes^{2}}_{1}(x\otimes y) -e_{2}\otimes\mu_{1}\circ\alpha^{\otimes^{2}}_{1}(x\otimes y)\cr
&=&\alpha_{1}(\mu_{1}(x\otimes y))\otimes e_{1} + e_{1}\otimes\alpha_{1}(\mu_{1}(x\otimes y)) - e_{2}\otimes\alpha_{1}(\mu_{1}(x\otimes y))\cr
&=&\Delta_{1}(\mu_{1}(x\otimes y));
\eeqs
\beqs
\Delta_{1}(\alpha_{1}(x))&=&\alpha^{2}_{1}(x)\otimes e_{1} + e_{1}\otimes\alpha^{2}_{1}(x) - e_{2}\otimes\alpha^{2}_{1}(x)\cr
&=&\alpha^{\otimes^{2}}_{1}(\alpha_{1}(x)\otimes e_{1} + e_{1}\otimes\alpha_{1}(x) - e_{2}\otimes\alpha_{1}(x))\cr
&=& \alpha^{\otimes^{2}}_{1}\circ\Delta_{1}(x).
\eeqs
Then, $\Delta$ is a homomorphism of the hom-associative algebras $(V, \mu, \alpha)$ and $(V\otimes V, \bullet, \alpha\otimes\alpha)$. 
 \end{enumerate}
Therefore, we can conclude that $\mathcal{K}_{1}(\mathcal{A}):=(\tilde{V}, \mu_{1}, \eta_{1}, \Delta_{1}, \varepsilon_{1}, \alpha_{1})$ is a unital hom-bialgebra.

$ \hfill \square $

 \begin{proposition}\label{k2}
 Let $\mathcal{A}=(V, \mu, \eta, \alpha)$ be a unital hom-associative algebra, (where $e_{2}:= \eta(1)$ is the unit). Let $\tilde{V}$ be the vector space spanned by $V$ and $e_{1}, \tilde{V}=span(V, e_{1})$. We have the unital hom-bialgebra $\mathcal{K}_{2}(\mathcal{A}):=(\tilde{V}, \mu_{2}, \eta_{2}, \Delta_{2}, \varepsilon_{2}, \alpha_{2}),$ where the multiplication $\mu_{2}$ is defined by:
 \beqs
&& \mu_{2}(e_{1}\otimes x)= \mu_{2}(x\otimes e_{1})= x \ \ \forall x\in \tilde{V},\cr
&& \mu_{2}(x\otimes y)= \mu(x\otimes y) \ \ \forall x, y\in V,
 \eeqs
 the unit $\eta_{2}$ is given by $\eta_{2}(1)= e_{1},$ while the comultiplication $\Delta_{2},$  the counit $\varepsilon_{2},$   and the linear map $\alpha_{2}$ are defined by:
 \beqs
 &&\Delta_{2}(e_{1})=e_{1}\otimes e_{1},\cr
 &&\Delta_{2}(e_{2})=e_{2}\otimes e_{1} + e_{1}\otimes e_{2} -e_{2}\otimes e_{2},\cr
 &&\Delta_{2}(x)= (e_{1}-e_{2})\otimes \alpha(x) + \alpha(x)\otimes (e_{1}- e_{2}) \ \ \forall x\in V\backslash \lbrace e_{2}\rbrace,
 \eeqs
 \beqs
 \varepsilon_{2}(e_{1})=1, \varepsilon_{2}(x)=0 \  \  \forall x\in V,
 \eeqs
 \beqs
 \alpha_{2}(e_{1})= e_{1}, \alpha_{2}(e_{2})= e_{2}, \alpha_{2}(x)=\alpha(x),\ \forall x\in V,
 \eeqs
respectively.
 \end{proposition}
 \textbf{Proof:}
 We have:
 \begin{enumerate}
 \item[$\diamond$]
 \beqs
\Delta_{2}(\alpha_{2}(x))&=&(e_{1}- e_{2})\otimes\alpha^{2}_{2}(x) + \alpha^{2}_{2}(x)\otimes(e_{1}-e_{2})\cr
&=&\alpha^{\otimes^{2}}_{2}((e_{1}- e_{2})\otimes\alpha_{2}(x) + \alpha_{2}(x)\otimes(e_{1}-e_{2}))\cr
&=&\alpha^{\otimes^{2}}_{2}(\Delta_{2}(x)).
 \eeqs
  \item[$\diamond$]
  \beqs
\Delta_{2}(x)\bullet \Delta_{2}(y)&=&((e_{1}- e_{2})\otimes\alpha_{2}(x) + \alpha_{2}(x)\otimes(e_{1}-e_{2}))\bullet ((e_{1}- e_{2})\otimes\alpha_{2}(y) \cr
&& + \alpha_{2}(y)\otimes(e_{1}-e_{2}))\cr
&=&\mu_{2}[(e_{1}- e_{2}); (e_{1}- e_{2})]\otimes\mu_{2}[\alpha_{2}(x); \alpha_{2}(y)]\cr &&+ \mu_{2}[(e_{1}- e_{2}); \alpha_{2}(y)]\otimes\mu_{2}[\alpha_{2}(x); (e_{1}-e_{2})]\cr
&& + \mu_{2}[\alpha_{2}(x); (e_{1}- e_{2})]\otimes\mu_{2}[(e_{1}- e_{2}); \alpha_{2}(y)]\cr
&& +\mu_{2}(\alpha_{2}(x); \alpha_{2}(y))\otimes\mu_{2}[(e_{1}- e_{2}); (e_{1}- e_{2})]\cr
&=&(e_{1}- e_{2})\otimes\mu_{2}(\alpha_{2}(x); \alpha_{2}(y)) + \mu_{2}(\alpha_{2}(x); \alpha_{2}(y))\otimes(e_{1}- e_{2})\cr
&=&(e_{1}- e_{2})\otimes\alpha_{2}(\mu_{2}(x, y)) + \alpha_{2}(\mu_{2}(x, y))\otimes(e_{1}- e_{2})\cr
&=&\Delta_{2}(\mu_{2}(x, y)).
  \eeqs
  \item[$\diamond$]
  \beqs
(\alpha_{2}\otimes\Delta_{2})\circ\Delta_{2}(x)&=&e_{1}\otimes e_{1}\otimes\alpha^{2}_{2}(x)- e_{1}\otimes e_{2}\otimes\alpha^{2}_{2}(x) -e_{2}\otimes e_{1}\otimes\alpha^{2}_{2}(x) \cr
&& + e_{2}\otimes e_{2}\otimes\alpha^{2}_{2}(x) + e_{1}\otimes\alpha^{2}_{2}(x)\otimes e_{1} -e_{1}\otimes\alpha^{2}_{2}(x)\otimes e_{2} -e_{2}\otimes\alpha^{2}_{2}(x)\otimes e_{2}\cr
&& + \alpha^{2}_{2}(x)\otimes e_{1}\otimes e_{1} -\alpha^{2}_{2}(x)\otimes e_{2}\otimes e_{1} -\alpha^{2}_{2}(x)\otimes e_{1}\otimes e_{2} + \alpha^{2}_{2}(x)\otimes e_{2}\otimes e_{2}\cr
&=& (\Delta_{2}\otimes\alpha_{2})\circ\Delta_{2}(x).
\eeqs
\item[$\diamond$] The condition (\ref{hom-counit condition}) is easily established.
 \end{enumerate}
Hence, $\mathcal{K}_{2}(\mathcal{A}):=(\tilde{V}, \mu_{2}, \eta_{2}, \Delta_{2}, \varepsilon_{2}, \alpha_{2})$ is a unital hom-bialgebra. 
 
 $ \hfill \square $
 \subsection{Construction of unital 2-hom-associative bialgebras}
 We construct  $(n+1)$-dimensional unital $2$-hom-associative bialgebras from  $n$-dimensional unital hom-associative algebras.
\begin{lemma}\label{lem1}
Let $\A=(V, \mu, \eta, \alpha)$ be any unital hom-associative algebra. The unital hom-bialgebra $\mathcal{K}_{1}(\A)=(\tilde{V}, \mu_{1}, \eta_{1}, \Delta_{1}, \varepsilon_{1}, \alpha_{1})$ is a unital infinitesimal hom-bialgebra.
\end{lemma}
\textbf{Proof:}
We know that $\mathcal{K}_{1}(\mathcal{A})$ is a unital hom-bialgebra. Then, we only have to show  the unital hom-infinitesimal condition. For all $x, y\in V,$ we have:
\beqs
&&(\mu\otimes\alpha)\circ(\alpha\otimes\Delta)(x\otimes y) + (\alpha\otimes\mu)\circ(\Delta\otimes\alpha)(x\otimes y) -\alpha^{2}(x)\otimes\alpha(y)\cr
&&=(\mu\otimes\alpha)[\alpha(x)\otimes(\alpha(y)\otimes e_{1} + e_{1}\otimes\alpha(y) -e_{2}\otimes\alpha(y))] + (\alpha\otimes\mu)[(\alpha(x)\otimes e_{1} + e_{1}\otimes\alpha(x)\cr
&& -e_{2}\otimes\alpha(x))\otimes\alpha(y)] -\alpha^{2}(x)\otimes\alpha(y)\cr
&&=(\mu\otimes\alpha)[\alpha(x)\otimes\alpha(y)\otimes e_{1} + \alpha(x)\otimes e_{1}\otimes\alpha(y) -\alpha(x)\otimes e_{2}\otimes\alpha(y)] + \cr
&&(\alpha\otimes\mu)[\alpha(x)\otimes e_{1}\otimes\alpha(y) + e_{1}\otimes\alpha(x)\otimes\alpha(y)  -e_{2}\otimes\alpha(x)\otimes\alpha(y)] -\alpha^{2}(x)\otimes\alpha(y)\cr
&&=\mu(\alpha(x), \alpha(y))\otimes e_{1} + \alpha(x)\otimes\alpha^{2}(y) - \alpha(x)\otimes\alpha^{2}(y) + \alpha^{2}(x)\otimes\alpha(y) + e_{1}\otimes\mu(\alpha(x), \alpha(y))\cr
&&-e_{2}\otimes\mu(\alpha(x), \alpha(y)) - \alpha^{2}(x)\otimes\alpha(y)\cr
&&=\mu(\alpha(x), \alpha(y))\otimes e_{1} + e_{1}\otimes\mu(\alpha(x), \alpha(y)) -e_{2}\otimes\mu(\alpha(x), \alpha(y))\cr
&&=\alpha(\mu(x, y))\otimes e_{1} + e_{1}\otimes \alpha(\mu(x, y)) - e_{2}\otimes\alpha(\mu(x, y))\cr
&&=\Delta(\mu(x, y)).
\eeqs
Hence, $\mathcal{K}_{1}(\A)$ is a unital infinitesimal hom-bialgebra.
 
 $ \hfill \square $
\begin{remark}
Let $\A=(V, \mu, \eta, \alpha)$ be any unital hom-associative algebra. The unital hom-bialgebra $\mathcal{K}_{2}(\A)$ is not a unital infinitesimal hom-bialgebra since the unital hom-infinitesimal condition is not satisfied.
\end{remark}
\begin{remark}
Let $\A_{2}=(V, \mu_{1}, \mu_{2}, \eta, \alpha)$ be a unital $2$-hom-associative algebra, then we have the same hom-coalgebra structure in the associated hom-bialgebra, (or  unital infinitesimal hom-bialgebra),  related to initial  unital hom-associative algebras $(V, \mu_{1},  \eta, \alpha)$ and $(V,  \mu_{2}, \eta, \alpha)$.
\end{remark}
\begin{proposition}
Let $\A=(V, \mu, \eta, \alpha)$ and $\A'=(V, \mu', \eta, \alpha)$ be  two unital hom-associative algebras over an $n$-dimensional vector space $V$. Let $\mathcal{K}_{1}(\A)=(\tilde{V}, \mu_{1}, \eta_{1}, \Delta_{1}, \varepsilon_{1}, \alpha_{1})$ and $\mathcal{K}_{1}(\A')=(\tilde{V}, \mu'_{1}, \eta_{1}, \Delta_{1}, \varepsilon_{1}, \alpha_{1})$ be the above defined associated hom-bialgebras. Then, we have that $\mathcal{B}_{1}= (\tilde{V}, \mu_{1}, \mu'_{1}, \eta_{1}, \Delta_{1}, \varepsilon_{1}, \alpha_{1})$ is a $(n + 1)$-dimensional unital 2-hom-associative bialgebra over the vector space $\tilde{V}= span(V, e_{1}),$ where $\eta_{1}(1)= e_{1}$.
\end{proposition}
\textbf{Proof}

From Lemma \ref{lem1}, $\mathcal{K}_{1}(\A')$ is a unital infinitesimal hom-bialgebra, and $\mathcal{K}_{1}(\A)$ is a unital hom-bialgebra, and hence $\mathcal{B}_{1}= (\tilde{V}, \mu_{1}, \mu'_{1}, \eta_{1}, \Delta_{1}, \varepsilon_{1}, \alpha_{1})$ is a unital 2-hom-asssociative bialgebra.

$ \hfill \square $
\begin{remark}
Let $\mathcal{B}=(V, \mu, \eta, \Delta, \varepsilon, \alpha)$ be a unital hom-bialgebra.  If the comultiplication satisfies the unital hom-infinitesimal condition, then $\mathcal{B}=(V, \mu, \mu,\eta, \Delta, \varepsilon, \alpha)$ is a unital $2$-associative-hom-bialgebra.
\end{remark}
\subsection{Construction of unital $2$-hom bialgebras}
\begin{proposition}
Let $V$ be an $n$-dimensional vector space over $\K$. Let $\A_{1}=(V, \mu_{1}, \eta_{1}, \alpha)$ and $\A_{2}=(V, \mu_{2}, \eta_{2}, \alpha)$ be two unital hom-associative algebras, and  $\mathcal{K}_{j}(\A_{i})=(\tilde{V}, \tilde{\mu}_{i}, \eta, \Delta_{i}, \varepsilon, \tilde{\alpha}), i, j=1, 2$  the above defined associated hom-bialgebras. Then, 
\beqs
\mathcal{B}_{1}=(\tilde{V}, \tilde{\mu}_{1}, \tilde{\mu}_{2}, \eta, \Delta_{1}, \Delta_{2}, \varepsilon, \tilde{\alpha}) \mbox{ and } \mathcal{B}_{2}=(\tilde{V}, \tilde{\mu}_{1}, \tilde{\mu}_{2}, \eta, \Delta^{cop}_{1}, \Delta_{2}, \varepsilon, \tilde{\alpha})
\eeqs
are two $(n+1)$-dimensional unital $2$-hom bialgebras on $\tilde{V}=span(V, e_{1}),$ where $\eta(1)= e_{1}$
\end{proposition}
\textbf{Proof:} Similarly to the 
 proofs of Theorems \ref{k1} and \ref{k2}, we establish by a straightforward computation that  $\mathcal{B}_{1}$ and $\mathcal{B}_{2}$ are unital 2-hom-bialgebras.

$ \hfill \square $

The next corollary  gives a unital $2$-$2$-hom-bialgebra  from  two unital hom-associative algebras.
\begin{corollary}
Let $V$ be an $n$-dimensional vector space over $\K$. Let $\A_{1}=(V, \mu_{1}, \eta_{1}, \alpha)$ and $\A_{2}=(V, \mu_{2}, \eta_{2}, \alpha)$ be two unital hom-associative algebras, and $\mathcal{K}_{1}(\A_{i}), i=1, 2$ be the above defined associated bialgebras.
Then, $\mathcal{B}_{1}=(\tilde{V}, \tilde{\mu}_{1}, \tilde{\mu}_{2}, \eta, \Delta_{1}, \Delta_{2}, \varepsilon, \tilde{\alpha})$ is a $(n+1)$-dimensional unital $2$-$2$-hom bialgebra on $\tilde{V}=span(V, e_{1}),$ where $\eta(1)= e_{1}.$
\end{corollary}
\section{Hom-left symmetric dialgebras}
In this section, we discuss  hom-associative dialgebras and  hom-left symmetric dialgebras generalizing the classical left-symmetric dialgebras in the   hom-algebra setting.
\subsection{Hom-associative dialgebra}
The associative dialgebras were introduced by Loday in \cite{[Loday2]}. The notion of dialgebra is a generalization
 of an associative algebra with two operations, which gives rise to a Leibniz algebra instead of a Lie algebra.
\begin{definition}
A hom-associative dialgebra is a quadruple $(D, \dashv, \vdash, \alpha)$ consisting of a linear space $D$, two associative products $\dashv, \vdash: D\times D\rightarrow D,$ and a homomorphism $\alpha: D\rightarrow D$ satisfying
\beq\label{D1}
\alpha(x)\dashv(y\dashv z)= \alpha(x)\dashv(y\vdash z)
\eeq
\beq\label{D2}
(x\vdash y)\dashv\alpha(z)= \alpha(x)\vdash(y\dashv z)
\eeq
\beq\label{D3}
(x\vdash y)\vdash \alpha(z)=(x\dashv y)\vdash\alpha(z).
\eeq 
\end{definition}
We recover the classical associative dialgebra when $\alpha= \id.$
\begin{example}
Let $\A$ be a hom-associative algebra, then the products $x\dashv y=xy= x\vdash y$ define a structure of hom-associative dialgebra on $\A$.
\end{example}
\begin{definition} 
We call differential hom-associative algebra the quadruple $(\A, \cdot, \alpha, d)$
 such that  $(\A, \cdot, \alpha)$ is a hom-associative algebra, $d(a\cdot b)= da\cdot b + a\cdot
  db, \forall a, b\in \A$, $d^{2}= 0,$ and $d\circ\alpha= \alpha\circ d$. 
\end{definition}
  \begin{proposition}
Let $(\A, \cdot, \alpha, d)$ be a differential hom-associative algebra. Consider the products    $\vdash$ and $\dashv$ on $\A$ given by $x\dashv y= \alpha(x)d\alpha(y)$ and $x\vdash y=\alpha(x)d\alpha(y)$. Then, $(\A, \dashv,
    \vdash, \alpha)$ is a hom-associative dialgebra.
\end{proposition}
\textbf{Proof:}
By hypothesis, $(\A, \cdot, \alpha, d)$ is a differential hom-associative algebra. Hence we have:
\begin{enumerate}
\item[$\bullet$]
$
\alpha(x)\dashv(y\dashv z)
= \alpha(x)\dashv(\alpha(y)d\alpha(z))=\alpha^{2}
(x)d\alpha(\alpha(y)d\alpha(z))=\alpha^{2}(x)d[\alpha^{2}(y)\alpha(d\alpha(z))]\\
=\alpha^{2}(x)d\alpha^{2}(y)d\alpha^{2}(z)= \alpha^{2}(x)d[\alpha(d\alpha(y))\alpha^{2}(z)]= \alpha^{2}
(x)\dashv\alpha[d\alpha(y)\alpha(z)]\\
= \alpha(x)\dashv(d\alpha(y)\alpha(z))
=\alpha(x)\dashv(y\vdash z);
$
\item[$\bullet$]
$
(x\vdash y)\dashv\alpha(z)
=(d\alpha(x)\alpha(y))\dashv\alpha(z)=d\alpha^{2}(x)\alpha^{2}(y)d\alpha^{2}(z)\\= \alpha(x)\vdash(\alpha(y)d\alpha(z))=\alpha(x)\vdash(y\dashv z);
$
\item[$\bullet$]
$
(x\vdash y)\vdash\alpha(z)
=(d\alpha(x)\alpha(y))\vdash\alpha(z)=d[d\alpha^{2}(x)\alpha^{2}(y)]\alpha^{2}(z)= d\alpha^{2}(x)d\alpha^{2}(y)\alpha^{2}(z)\\=d(\alpha^{2}(x)d\alpha^{2}(y))\alpha^{2}(z)= (\alpha(x)d\alpha(y))\vdash\alpha(z)=(x\dashv y)\vdash\alpha(z).
$
\end{enumerate} 
Therefore, $(A, \dashv, \vdash, \alpha)$ is a hom-associative dialgebra.
$ \hfill \square $
\begin{definition}
Let $(D, \dashv, \vdash, \alpha)$ and $(D', \dashv', \vdash', \alpha')$ be two hom-associative dialgebras. A linear map $f: D\rightarrow D'$ is a hom-associative dialgebras morphism if 
\beqs
\dashv'\circ(f\otimes f)= f\circ \dashv,\ \vdash'\circ(f\otimes f)= f\circ \vdash \mbox{ and }  f\circ\alpha= \alpha'\circ f.
\eeqs
\end{definition}
Now, we show that we can construct hom-associative dialgebras starting from a classical associative dialgebra and an algebra endomorphism. We extend then the construction by the composition introduced by D. Yau in \cite{[yau2]} for Lie and associative algebras.
\begin{theorem}
Let $(D, \dashv, \vdash)$ be an associative dialgebra and $\alpha: D\rightarrow D$ be an
 associative dialgebra endomorphism. Then $D_{\alpha}=(D, \dashv_{\alpha}, \vdash_{\alpha},
  \alpha)$, where $\dashv_{\alpha}= \alpha\circ\dashv$ and $\vdash_{\alpha}= \alpha\circ\vdash$, is
   a hom-associative dialgebra. Moreover, suppose that $(D', \dashv', \vdash')$ is another
    associative dialgebra and $\alpha': D'\rightarrow D'$ is an associative dialgebra endomorphism.
     If $f: D\rightarrow D'$ is an associative dialgebra morphism that satisfies $f\circ\alpha=
      \alpha'\circ f,$ then $f: D_{\alpha}\rightarrow D'_{\alpha'}$ is a morphism of hom-associative dialgebras.
\end{theorem}
\textbf{Proof:}
We have:
\begin{enumerate}
\item[$\bullet$] 
$
\alpha(x)\dashv_{\alpha}(y\dashv_{\alpha}z)
= \alpha(\alpha(x)\dashv (y\dashv_{\alpha}z))= \alpha(\alpha(x)\dashv \alpha(y\dashv z))=\alpha^{2}(x\dashv (y\dashv z))\\
= \alpha^{2}(x\dashv (y\vdash z))= \alpha(\alpha(x)\dashv(\alpha(y\vdash z)))= \alpha(\alpha(x)\dashv(y\vdash_{\alpha} z))
= \alpha(x)\dashv_{\alpha}(y\vdash_{\alpha} z);
$
\item[$\bullet$]
$
(x\vdash_{\alpha} y)\dashv_{\alpha}\alpha(z)=\alpha^{2}((x\vdash y)\dashv z)= \alpha^{2}(x\vdash (y\dashv z))= \alpha(x)\vdash_{\alpha}(y\dashv_{\alpha} z);
$
\item[$\bullet$]
$
(x\vdash_{\alpha} y)\vdash_{\alpha} \alpha(z)= \alpha^{2}((x\vdash y)\vdash z) = \alpha^{2}((x\dashv y)\vdash z) =(x\dashv_{\alpha} y)\vdash_{\alpha}\alpha(z);
$

\item[$\bullet$]
$f\circ\dashv_{\alpha}
= f\circ(\alpha\circ \dashv)= (f\circ\alpha)\circ\dashv=(\alpha'\circ f)\circ\dashv= \alpha'\circ(f\circ\dashv)=\alpha'\circ(\dashv'\circ(f\otimes f))\\= (\alpha'\circ \dashv')\circ(f\otimes f)= \dashv'_{\alpha'}\circ(f\otimes f),$ and we also obtain that $f\circ\vdash_{\alpha}= \vdash'_{\alpha'}\circ(f\otimes f).$
\end{enumerate}
Hence, we can conclude that $D_{\alpha}$ is a hom-associative dialgebra and $f$ a morphism of hom-associative dialgebras. 

$ \hfill \square $
\subsection{Hom-Leibniz algebra}
\begin{definition}
A hom-Leibniz algebra is a triple $(L, [.,.] \alpha)$ consisting of a linear space $L$, a bilinear product $[.,.]: L\times L\rightarrow L,$ and a  homomorphism $\alpha: L\rightarrow L$ satisfying
\beq\label{identity hom-leibniz algebra}
[[x, y], \alpha(z)]= [[x, z], \alpha(y)] + [\alpha(x), [y, z]].
\eeq
If the bracket is skew-symmetric, then $L$ is a hom-Lie algebra. Hence, the  hom-Lie algebras are particular cases of hom-Leibniz algebras.
\end{definition}
\begin{proposition}
Let $(\A, \cdot, \alpha, d)$ be a differential hom-associative algebra. Defining the bracket on $\A$  by 
\beqs
[x, y]:= \alpha(x)\cdot d\alpha(y)-d\alpha(y)\cdot \alpha(x), 
\eeqs
then the vector space $\A$ equipped with this bracket is a hom-Leibniz algebra. 
\end{proposition}
\textbf{Proof:}
By direct computation, we obtain:
\begin{enumerate}
\item[$\bullet$]
$
[\alpha(x), [y, z]]= \alpha(x)d\alpha(y)d\alpha(z)- d\alpha(y)d\alpha(z)\alpha(x)
 -\alpha(x)d\alpha(z)d\alpha(y) + d\alpha(z)d\alpha(y)\alpha(x);
$
\item[$\bullet$]
$
[[x, y], \alpha(z)]= \alpha(x)d\alpha(y)d\alpha(z)-
 d\alpha(z)\alpha(x)d\alpha(y)-
d\alpha(y)\alpha(x)d\alpha(z) + d\alpha(z)d\alpha(y)\alpha(x);
$
\item[$\bullet$]
$
[[x, z], \alpha(y)]= \alpha(x)d\alpha(z)d\alpha(y)-d\alpha(y)\alpha(x)d\alpha(z)-d\alpha(z)\alpha(x)d\alpha(y)+ d\alpha(y)d\alpha(z)\alpha(x).
$
\end{enumerate}
Then, $[\alpha(x), [y, z]]= [[x, y], \alpha(z)] - [[x, z], \alpha(y)]$. Hence, the pair $(A, [.,.])$
is a hom-Leibniz algebra.

$ \hfill \square $

A morphism of hom-Leibniz algebras is a linear map of the underlying vector space that commutes with the twisting maps and the multiplications.
\begin{theorem}
Let $(L, [.,.])$ be a Leibniz algebra, and $\alpha: L\rightarrow L$ be a Leibniz algebra endomorphism. Then $L_{\alpha}=(L, [.,.]_{\alpha}, \alpha)$ is a hom-Leibniz algebra. Moreover, suppose that $(L', [.,.]')$ is another Leibniz algebra, and $\alpha': L'\rightarrow L'$ a Leibniz algebra endomorphism. If $f: L\rightarrow L'$ is a Leibniz algebra morphism  satisfying $f\circ\alpha= \alpha'\circ f,$ then $L_{\alpha}\rightarrow L'_{\alpha'}$ is a morphism of Leibniz algebras.
\end{theorem}
\textbf{Proof:}
Since
\begin{enumerate}
\item[$\bullet$]
$
[[x, y]_{\alpha}, \alpha(z)]_{\alpha}
= \alpha([\alpha([x, y]), \alpha(z)])= \alpha^{2}([[x, y], z])= \alpha^{2}([[x, z], y] + [x, [y, z]])\\ = [[x, z]_{\alpha}, \alpha(y)]_{\alpha} + [\alpha(x), [y, z]_{\alpha}]_{\alpha}
$ and
\item[$\bullet$]
$
f\circ[.,.]_{\alpha}
= f\circ(\alpha\circ [.,.])= (f\circ \alpha)\circ [.,.]= (\alpha'\circ f)\circ[.,.]=\alpha'\circ(f\circ [.,.])\\ =\alpha'\circ([.,.]'\circ (f\otimes f))= (\alpha'\circ[.,.]')\circ(f\otimes f)= [.,.]'_{\alpha'}\circ (f\otimes f),
$
\end{enumerate}
then we have the results.

$ \hfill \square $

\begin{theorem}
Let $(D, \dashv, \vdash, \alpha)$ be a hom-associative dialgebra. Let $[.,.]: D\otimes D\rightarrow 
D$ be a linear map defined for $x, y\in D$ by 
\beqs
[x,y]= x\dashv y - x\vdash y.
\eeqs
Then $(D, [.,.], \alpha)$ is a hom-Leibniz algebra.
\end{theorem}
\textbf{Proof:}
$(D, \dashv, \vdash, \alpha)$ is a hom-associative dialgebra, then we have
\beqs
\alpha(y)\vdash(z\vdash x)=(y\vdash z)\vdash\alpha(x)= (y\dashv z)\vdash\alpha(x);\cr
(x\dashv z)\dashv\alpha(y)= \alpha(x)\dashv(z\dashv y)=\alpha(x)\dashv(z\vdash y).
\eeqs
Therefore, by direct computation, we obtain the identity \ref{identity hom-leibniz algebra}.

$ \hfill \square $
\subsection{Hom-left symmetric dialgebras}
Now, we generalize the notion of left symmetric dialgebra introduced by R. Felipe,  twisting the identities by a linear map, as well as some theorems established in \cite{[Felipe2]}.
\begin{definition}
A left hom-preLie algebra is a triple $(\A, ., \alpha)$ consisting of a vector space $\A$, a bilinear map $.: \A\otimes \A\rightarrow \A,$ and a homomorphism $\alpha$ satisfying 
\beq
\alpha(x).(y.z) - (x.y).\alpha(z)= \alpha(y).(x.z) - (y.x).\alpha(z).
\eeq
\end{definition}
\begin{remark}
Any hom-associative algebra is a left hom-preLie algebra.
\end{remark}
A morphism of left hom-preLie algebras is a linear map of the underlying vector space that commutes with the twisting maps and the multiplications.
\begin{definition}
Let $S$ be a vector space over a field $K$. Let us assume that $S$ is equipped with two biliniear products $\dashv, \vdash:S\otimes S\rightarrow S,$ and a  homomorphism $\alpha: S\rightarrow S$  satisfying the identities:
\beq\label{hlsdeq1}
\alpha(x)\dashv(y\dashv z)&=& \alpha(x)\dashv(y\vdash z),
\eeq
\beq\label{hlsdeq2}
(x\vdash y)\vdash\alpha(z)&=&(x\dashv y)\vdash\alpha(z),
\eeq
\beq\label{hlsdeq3}
\alpha(x)\dashv(y\dashv z)- (x\dashv y)\dashv\alpha(z)&=& \alpha(y)\vdash(x\dashv z) - (y\vdash x)\dashv\alpha(z),
\eeq
\beq\label{hlsdeq4}
\alpha(x)\vdash(y\vdash z)- (x\vdash y)\vdash\alpha(z)&=&\alpha(y)\vdash(x\vdash z)- (y\vdash x)\vdash\alpha(z).
\eeq
Then, we say that $S$ is a hom-left symmetric dialgebra (HLSDA), or left disymmetric hom-algebra.
\end{definition}
\begin{example}
Any hom-associative algebra $(\A, \cdot, \alpha)$ is a hom-left symmetric dialgebra with $\vdash=\dashv=\cdot$.
\end{example}
The definition of a hom-left symmetric dialgebras morphism is similar to that of a hom-associative dialgebras morphism. 
\begin{remark}
We can  construct a hom-left symmetric dialgebra by the composition method  from a classical left-symmetric dialgebra $(D, \dashv, \vdash)$ and an algebra endomrophism $\alpha,$ by considering $(D, \dashv_{\alpha}, \vdash_{\alpha}, \alpha),$ where $x\dashv_{\alpha} y= \alpha(x\dashv y)$ and $x\vdash_{\alpha} y= \alpha(x\vdash y)$.
\end{remark}
Let us denoted by $\mathcal{HS}$ the set of all hom-left symmetric dialgebras, and $\mathcal{HD}$ the set of all hom-associative dialgebras.
\begin{proposition}
All hom-associative dialgebra is a hom-left symmetric dialgebra.

 Then,$\mathcal{HD}\subseteq\mathcal{HS}.$
\end{proposition}
\textbf{Proof:}
Let $(D, \dashv, \vdash, \alpha)$ be a hom-associative dialgebra, then the Eqs. (\ref{hlsdeq1}) and (\ref{hlsdeq2}) are satisfied. Since the products $\dashv$ and $\vdash$ are associative with the condition (\ref{D2}), hence the Eqs. (\ref{hlsdeq3}) and (\ref{hlsdeq4}) are established.
   
$ \hfill \square $
\begin{remark}
Any hom-left symmetric algebra is a hom-left symmetric dialgebra in which $\vdash= \dashv.$ A nonassociative hom-left symmetric algebra is not a hom-left symmetric dialgebra. Hence, we have $\mathcal{HD}\neq\mathcal{HS}.$
\end{remark}
\begin{proposition}
A hom-left symmetric dialgebra $S$ is a hom-associative dialgebra if and only if both products of $S$ are hom-associative.
\end{proposition}
\textbf{Proof:}
Let $(S, \dashv, \vdash, \alpha)$ be a hom-left symmetric dialgebra. If $S$ is a hom-associative dialgebra, then the products $\dashv$ and $\vdash$ defined on $S$ are hom-associative. Conversely, suppose that the products $\dashv$ and $\vdash$ are hom-associative. Since $S$ has a hom-left symmetric dialgebra structure, then from the Eq. (\ref{hlsdeq3}), the Eq. (\ref{D2}) is satisfied. Hence $S$ is a hom-associative dialgebra.

$ \hfill \square $

\begin{theorem}
Let $(S, \vdash, \dashv, \alpha)$ be a hom-left symmetric dialgebra. Then, the commutator given by $[x, y]=x\dashv y- y\vdash x$ defines a structure of hom-leibniz algebra on $S.$ In others words, $(S, [.,.], \alpha)$ is a hom-Leibniz algebra.
\end{theorem}
\textbf{Proof:}
We have:
\begin{enumerate}
\item[$\bullet$]
$
[[x, y], \alpha(z)]= (x\dashv y)\dashv\alpha(z) - \alpha(z)\vdash(x\dashv y) - (y\vdash x)\dashv\alpha(z) + \alpha(z)\vdash(y\vdash x);
$
\item[$\bullet$]
$
[[x, z], \alpha(y)]=(x\dashv z)\dashv\alpha(y)-\alpha(y)\vdash(x\dashv z) + (z\vdash x)\dashv\alpha(y) - \alpha(y)\vdash(z\vdash x);
$
\item[$\bullet$]
$
[\alpha(x),[y, z]]=\alpha(x)\dashv(y\dashv z) - (y\dashv z)\vdash\alpha(x)- \alpha(x)\dashv(z\vdash y) + (z\vdash y)\vdash\alpha(x).
$
\end{enumerate}
 From the Eqs.(\ref{hlsdeq3}) and (\ref{hlsdeq4}), we obtain the condition \ref{identity hom-leibniz algebra}.

$\hfill \square $


\begin{definition}
Let $(L, [.,.], \alpha)$ be a hom-Leibniz algebra. The pair of bilinear mappings $\bigtriangledown_{1}, \bigtriangledown_{2}: L\times L\rightarrow L$ is called an affine hom-Leibniz structure obeying the relations:
\beq
\bigtriangledown_{2}(x, y) - \bigtriangledown_{1}(y, x)= [x, y],
\eeq
\beq \label{aff.leib.eq1}
\bigtriangledown_{1}(\bigtriangledown_{1}(x, y), \alpha(z))= \bigtriangledown_{1}(\bigtriangledown_{2}(x, y), \alpha(z)); \bigtriangledown_{2}(\alpha(x), \bigtriangledown_{2}(y, z))= \bigtriangledown_{2}(\alpha(x), \bigtriangledown_{1}(y, z))
\eeq
\beq \label{aff.leib.eq2}
\bigtriangledown_{2}(\alpha(x), \bigtriangledown_{2}(y; z)) - \bigtriangledown_{1}(\alpha(y), \bigtriangledown_{2}(x, z))= \bigtriangledown_{2}([x, y], \alpha(z))
\eeq
and
\beq\label{aff.leib.eq3}
\bigtriangledown_{1}(\alpha(x), \bigtriangledown_{1}(y, z)) - \bigtriangledown_{1}(\alpha(y), \bigtriangledown_{1}(x, z))= \bigtriangledown_{1}([x, y], \alpha(z))
\eeq
for all $x, y, z\in L.$
\end{definition}
\begin{theorem}
Let $(L, [.,.], \alpha)$ be a hom-Leibniz algebra and let $\bigtriangledown_{1}, \bigtriangledown_{2}$ define an affine hom-Leibniz structure. Then, $L$ is a hom-left symmetric dialgebra with $\vdash$ and $\dashv$ defined as 
\beq
x\vdash y=\bigtriangledown_{1}(x, y);\ x\dashv y=\bigtriangledown_{2}(x, y). 
\eeq
\end{theorem}
\textbf{Proof:}
(\ref{aff.leib.eq1}) implies (\ref{hlsdeq1}) and (\ref{hlsdeq2}). Then,  (\ref{hlsdeq3})
and (\ref{hlsdeq4}) follow from (\ref{aff.leib.eq2}) and (\ref{aff.leib.eq3}), respectively.

$\hfill \square $
\section{Concluding Remarks}
In this work, from the hom-counital  and  unital infinitesimal hom-bialgebra conditions, and following  Kaplansky's constructions based on  unital hom-associative algebras, we have built   unital 2-hom-associative bialgebras,  unital 2-hom-bialgebras, and unital 2-2-hom-bialgebras, and discussed their respective properties. 
Finally, we have defined and characterized the hom-left symmetric dialgebras  generalizing the
ordinary left symmetric dialgebras.
 \section*{Aknowledgement}
 The
 ICMPA is grateful to the  partnership with the Daniel Iagolnitzer Foundation (DIF), France, which supports the development of mathematical physics in Africa.

\end{document}